\documentclass[12pt,a4paper]{article}
\usepackage{latexsym,amsmath}
\usepackage{amssymb}
\usepackage{mathrsfs}
\usepackage{tocbibind}
\usepackage{cite}
\usepackage{hyperref}
\usepackage{amsmath}
\usepackage{tabularx}
\usepackage{array}
\usepackage[numbers,sort&compress]{natbib}

\nonstopmode \numberwithin{equation}{section}

\begin{document}
\begin{center}
\textbf{\large\bf Truncated Gauss hypergeometric series and its application in digamma function} \\
\end{center}

\begin{center}
{\bf M.I. Qureshi, Saima Jabee and M. Shadab{$^*$}} \\
\vskip.2cm
Department of Applied Sciences and Humanities,\\
Faculty of Engineering and Technology,Jamia Millia Islamia\\
(A Central University), New Delhi-110025,India.\\
{E-mails: miqureshi\_delhi@yahoo.co.in, saimajabee007@gmail.com,\\
 and shadabmohd786@gmail.com.}\\
{$^*$}{\bf Corresponding author}
\end{center}

{\bf Abstract:} In the last decades, the theory of digamma function has been developed with a high impact of interest by many authors. Here, we established some interesting results for digamma function, and also we have computed the values of digamma function for positive integers, using the concept of hypergeometric series. An attempt has been made to present some summation theorems for Clausen's hypergeometric function. Results presented here are potentially useful in the further study of digamma function.


\renewcommand{\thefootnote}{}


\footnote{Euler-Mascheroni constant; Gamma function. }

\section{Introduction}

A natural property of digamma function to be used as application in the theory of beta distributions-probability models for the domain [0,1]. It is mainly used in the theory of special functions with a wide range of the applications. Digamma functions are directly connected with many special functions such as Riemann's zeta function and Clausen's function etc.\\
Authors, who have participated in the theory of development of digamma function with respect to properties \cite{Qi, Gauss, Jensen, Lehmer, Mahler}, inequalities \cite{Alzer, Batir, Clark}, monotonicity \cite{Qi1, Qi2, Qi3, Qi4}, series \cite{Bor, de, Lewin, Wu, Gosper, Grossman}, and fractional calculus \cite{Al-saqabi, Sri1, Sri2}.
\vskip.2cm
The Gamma function, $\Gamma{(z)}$, was introduced by Leonard Euler as a generalization of the factorial function on the sets, $\mathbb{R}$ of all real numbers, and $\mathbb{C}$ of all complex numbers. It (or, Euler's integral of the second kind) is defined by

\begin{eqnarray}\label{eq(1.1)}
\Gamma{(z)}&=&\int_{0}^{\infty}\exp{(-t)}t^{z-1}dt,~~~~~~~~ \Re(z)>0\nonumber\\
&=&\lim_{n\rightarrow\infty}\int_{0}^{n}{\left(1-\frac{t}{n}\right)^n}t^{z-1}dt.
\end{eqnarray}
\vskip.2cm
In 1856, Karl Weierstrass gave a novel definition of gamma function
\begin{eqnarray}\label{eq(1.2)}
\frac{1}{\Gamma{(z)}}&=&z \exp{(\gamma z)} \prod_{n=1}^{\infty}\left[\left(1+\frac{z}{n}\right)\exp{\left(-\frac{z}{n}\right)}\right],
\end{eqnarray}
where $\gamma= 0.577215664901532860606512090082402431042\dots$ is called Euler-Mascheroni constant, and
$\frac{1}{\Gamma{(z)}}$ is an entire function of $z$, and
\begin{eqnarray}
\gamma&=&\lim_{n\rightarrow\infty}\left(1+\frac{1}{2}+\frac{1}{3}+.....+\frac{1}{n}-\ln {n}\right)\nonumber.
\end{eqnarray}

The function $\psi(z)$ is the logarithmic derivative of the gamma function or digamma function or Psi-function, given by\\

\begin{eqnarray}
\psi(z)=\frac{d}{dz}\{\ln{\Gamma(z)}\}=\frac{\Gamma ^{\prime}(z)}{\Gamma(z)},
\end{eqnarray}

\begin{eqnarray}
\ln{\Gamma(z)}=\int_{1}^{z}\psi(\zeta)d\zeta.
\end{eqnarray}

The widely-used Pochhammer symbol $(\lambda)_{\nu}$ ~$(\lambda, \nu \in\mathbb{C})$ is defined by
\begin{equation}
\left(\lambda\right)_{\nu}:=\frac{\Gamma\left(\lambda+\nu\right)}{\Gamma\left(\lambda\right)}=\begin{cases}
\begin{array}{c}
1\\
~\\
\lambda\left(\lambda+1\right)\ldots\left(\lambda+n-1\right)
\end{array} & \begin{array}{c}
\left(\nu=0;\lambda\in\mathbb{C}\setminus\left\{ 0\right\} \right)\\
~\\
\left(\nu=n\in\mathbb{N};\lambda\in\mathbb{C}\right) ,
\end{array}\end{cases}
\end{equation}
it is being understood $conventionally$ that $\left(0\right)_{0}=1$ and assumed $tacitly$ that the $\Gamma$ quotient exists.\\

\vskip.2cm

The {\it generalized hypergeometric function} ${_p}F_q$, is defined by

\begin{eqnarray}\label{eq(1.10)}
{_p}F_{q} \left( \begin{array}{r}(a_p); \\ (b_q); \end{array} z \right) = \sum_{m=0}^{\infty} \frac{[(a_p)]_m}{[(b_q)]_m} \frac{z^m}{m!}
\end{eqnarray}

\begin{itemize}
  \item p and q are positive integers or zero,
  \item $z$ is a complex variable,
  \item $(a_p)$ designates the set $\{a_1,a_2, . . . , a_p\}$,
  \item the numerator parameters $a_1, . . . , a_p \in\mathbb{C}$ and the denominator parameters \\$b_1, . . . , b_q \in\mathbb{C}\setminus\mathbb{Z}_{0}^{-}$ ,
  \item $[(a_r)]_k = \displaystyle\prod_{i=1}^{r} (a_i)_k$. By convention, a product over the empty set is 1.
  \end{itemize}

Thus, if a numerator parameter is a negative integer or zero, the $_pF_q$ series terminates, then we are led to a generalized hypergeometric polynomial.
\vskip.2cm

\section{Truncated Gauss series and its application in digamma function}

In 1931-32, W. N. Bailey (see, \cite[p. 40, last eq.]{Bailey2} and \cite[p. 34, last eq.]{Bailey1} ) derived a formula for truncated Gauss hypergeometric series in terms of Clausen's hypergeometric series as follows

\begin{eqnarray}\label{eq(2.6)}
{_2F_1} \left[ \begin{array}{r}a,b; \\ f; \end{array} 1 \right]_n  &=& \text{ Sum of first (n+1)-terms of series }{_2F_1} \left[ \begin{array}{r}a,b; \\ f; \end{array} 1 \right] \nonumber\\
&=& \sum_{k=0}^{n}\frac{(a)_k (b)_k }{(f)_k  k!} \nonumber\\
&=&\frac{\Gamma{(a+n+1)}\,\Gamma{(b+n+1)}}{\Gamma{(n+1)}\,\Gamma{(a+b+n+1)}}\,\,{_3}F_{2} \left[ \begin{array}{r}a,b,f+n; \\ f, a+b+n+1; \end{array} 1 \right],
\end{eqnarray}
where $ f \ge a+b $.
\vskip.1cm
It can be written as
\begin{eqnarray}\label{eq(2.2)}
{_3}F_{2} \left[ \begin{array}{r}a,b,f+n; \\ f, a+b+n+1; \end{array} 1 \right]=\frac{  \Gamma{(n+1)}\,\Gamma{(a+b+n+1)}}{\Gamma{(a+n+1)}\,\Gamma{(b+n+1)} } \sum_{k=0}^{n}\frac{(a)_k (b)_k }{(f)_k  k!},
\end{eqnarray}
where $ f \ge a+b $.

Putting $f = a+b+1$ in above equation \eqref{eq(2.2)}, we get

\begin{eqnarray}
{_2F_1} \left[ \begin{array}{r}a,b; \\a+b+1; \end{array} 1 \right]_n &=& \frac{\Gamma{(a+n+1)}\,\Gamma{(b+n+1)}}{\Gamma{(n+1)}\,\Gamma{(a+b+n+1)}}\,\,{_3}F_{2} \left[ \begin{array}{r}a,b,a+b+n+1; \\ a+b+1, a+b+n+1; \end{array} 1 \right]\nonumber\\
&=& \frac{\Gamma{(a+n+1)}\,\Gamma{(b+n+1)}}{\Gamma{(n+1)}\,\Gamma{(a+b+n+1)}}\,\,{_2}F_{1} \left[ \begin{array}{r}a,b; \\ a+b+1; \end{array} 1 \right]\nonumber\\
 \end{eqnarray}

Using Gauss summation theorem, we get
\begin{eqnarray}
 {_2F_1} \left[ \begin{array}{r}a,b; \\a+b+1; \end{array} 1 \right]_n = \frac{(a+1)_n \, (b+1)_n }{(a+b+1)_n \, n!}
 \end{eqnarray}

In \eqref{eq(2.2)}, put $a=1, b=1, f=2$ and $n=m-1$, we get

\begin{eqnarray}\label{eq(2.5)}
{_3F_2} \left[ \begin{array}{r}1,1,m+1 ;\\2,m+2; \end{array} 1 \right]=\left(\frac{m+1}{m} \right) \sum_{k=0}^{m-1} \frac{(1)_k}{(2)_k}=\left(\frac{m+1}{m} \right)\sum_{k=1}^{m} \frac{1}{k}, \qquad m =1,2,3,\dots,
\end{eqnarray}

it is the simplest way to calculate the Clausen's ${_3F_2}$-series for positive integers.\\

Now we calculate the values of Clausen's series ${_3F_2}$ by applying the result \eqref{eq(2.5)}.\\

In \eqref{eq(2.5)}, put $m=1$, we get
\begin{eqnarray}\label{eq(2.6)}
{_3F_2} \left[ \begin{array}{r}1,1,2 ;\\2,3; \end{array} 1 \right]=\left(\frac{2}{1} \right) \sum_{k=1}^{1} \frac{1}{k}=2.
\end{eqnarray}

In \eqref{eq(2.5)}, put $m=2$, we get
\begin{eqnarray}\label{eq(2.7)}
{_3F_2} \left[ \begin{array}{r}1,1,3 ;\\2,4; \end{array} 1 \right]=\left(\frac{3}{2} \right) \sum_{k=1}^{2} \frac{1}{k}=\frac{9}{4}.
\end{eqnarray}

In \eqref{eq(2.5)}, put $m=3$, we get
\begin{eqnarray}\label{eq(2.8)}
{_3F_2} \left[ \begin{array}{r}1,1,4 ;\\2,5; \end{array} 1 \right]=\left(\frac{4}{3} \right) \sum_{k=1}^{3} \frac{1}{k}=\frac{22}{9}.
\end{eqnarray}

Similarly for higher values of $m$ ($m=4,5,\dots,51$) in \eqref{eq(2.5)}, the values of ${_3F_2}$, are calculated using wolfram mathematica, and arranged in Table1, Table2 and Table3.

\vskip.2cm

Now, we establish an interesting formula for the computation of digamma function using Clausen's hypergeometric function given by \eqref{eq(2.5)}.\\

Let us recall the Weierstrass definition of gamma function \eqref{eq(1.2)}

\begin{eqnarray}
\frac{1}{\Gamma{(z)}}&=&z \exp{(\gamma z)} \prod_{n=1}^{\infty}\left[\left(1+\frac{z}{n}\right)\exp{\left(-\frac{z}{n}\right)}\right]\nonumber
\end{eqnarray}
Applying $\ln$ on both sides with respect to the base 'e'
\begin{eqnarray}
\ln{1}-\ln{\Gamma{(z)}}&=&\ln{z}+\gamma z\ln{e}+\sum_{n=1}^{\infty}\ln\left\{\left(1+\frac{z}{n}\right)\exp{\left(-\frac{z}{n}\right)}\right\}
\end{eqnarray}
\text{or,}
\begin{eqnarray}
-\ln{\Gamma{(z)}}&=&\ln{z}+\gamma z+\sum_{n=1}^{\infty}\left\{\ln\left(1+\frac{z}{n}\right)-\frac{z}{n}\right\}
\end{eqnarray}
Differentiating with respect to 'z', we get
\begin{eqnarray}\label{eq(2.11)}
\psi(z)=\frac{\Gamma{(z)^{'}}}{\Gamma{(z)}}&=&-\frac{1}{z}-\gamma -\sum_{n=1}^{\infty}\left\{\frac{1}{n+z}-\frac{1}{n}\right\}\nonumber\\
\psi(z)=\frac{\Gamma{(z)^{'}}}{\Gamma{(z)}}&=&-\frac{1}{z}-\gamma +\sum_{n=0}^{\infty}\left\{\frac{z}{(n+1)(n+z+1)}\right\}
\end{eqnarray}
Now writing R.H.S. of \eqref{eq(2.11)} in hypergeometric notation using \eqref{eq(1.10)}, we get
\begin{eqnarray}\label{eq(2.12)}
\psi(z)=\frac{\Gamma{(z)^{'}}}{\Gamma{(z)}}&=&-\frac{1}{z}-\gamma +\frac{z}{1+z}{_3}F_{2}\left[\begin{array}{r} 1, 1, 1+z;\\
~\\
2, 2+z ;\end{array}\  1\right],
\end{eqnarray}
 where $z\ne 0,-1,-2,-3,\dots$, and $\psi(z)$ denotes the digamma function.\\

For $z\in \mathbb{N}$ in \eqref{eq(2.12)}, and using the values of $_3F_2$ arranged in Tables \eqref{table-1.1},\eqref{table-1.2},\eqref{table-1.3}, we can find $\psi(1), \psi(2), \psi(3), \dots, \psi(52)$ (see Table \eqref{table-1.4}).

\newpage

\section{Some new summation theorems}

\begin{table}[h!]
\centering
\caption{( Clausen's summation theorem )}
\vskip.2cm
\label{table-1.1}
\begin{tabular}{|c  |c  |c  |c|c|c|c|c|c|}
\hline S.No. & Clasuen's series & S.No. & Clasuen's series \\ [5pt]\hline
1 & ${_3}F_{2}\left(\begin{array}{r} 1, 1, 3;\\
~\\
2, 4 ;\end{array}\  1\right) =\frac{9}{4}$ & 11 & ${_3}F_{2}\left(\begin{array}{r} 1, 1, 13;\\
~\\
2, 14 ;\end{array}\  1\right) = \frac{1118273}{332640}$     \\  [30pt]\hline

2 & ${_3}F_{2}\left(\begin{array}{r} 1, 1, 4;\\
~\\
2, 5 ;\end{array}\  1\right) =\frac{22}{9}$ &12 & ${_3}F_{2}\left(\begin{array}{r} 1, 1, 14;\\
~\\
2, 15;\end{array}\  1\right) = \frac{1145993}{334620}$     \\ [30pt]\hline

3 &  ${_3}F_{2}\left(\begin{array}{r} 1, 1, 5;\\
~\\
2, 6 ;\end{array}\  1\right) = \frac{125}{48}$ &13 & ${_3}F_{2}\left(\begin{array}{r} 1, 1, 15;\\
~\\
2, 16 ;\end{array}\  1\right) = \frac{1171733}{336336}$     \\ [30pt]\hline

4 & ${_3}F_{2}\left(\begin{array}{r} 1, 1, 6;\\
~\\
2, 7 ;\end{array}\  1\right) = \frac{137}{50}$ &14 & ${_3}F_{2}\left(\begin{array}{r} 1, 1, 16;\\
~\\
2, 17 ;\end{array}\  1\right) = \frac{2391514}{675675}$   \\ [30pt]\hline

5 &  ${_3}F_{2}\left(\begin{array}{r} 1, 1, 7;\\
~\\
2, 8 ;\end{array}\  1\right) = \frac{343}{120}$  &15 & ${_3}F_{2}\left(\begin{array}{r} 1, 1, 17;\\
~\\
2, 18 ;\end{array}\  1\right) = \frac{41421503}{11531520}$  \\ [30pt]\hline

6 &  ${_3}F_{2}\left(\begin{array}{r} 1, 1, 8;\\
~\\
2, 9 ;\end{array}\  1\right) = \frac{726}{245}$  &16 & ${_3}F_{2}\left(\begin{array}{r} 1, 1, 18;\\
~\\
2, 19 ;\end{array}\  1\right) = \frac{42142223}{11571560}$  \\ [30pt]\hline

7 &  ${_3}F_{2}\left(\begin{array}{r} 1, 1, 9;\\
~\\
2, 10 ;\end{array}\  1\right) = \frac{6849}{2240}$  &17 & ${_3}F_{2}\left(\begin{array}{r} 1, 1, 19;\\
~\\
2, 20 ;\end{array}\  1\right) = \frac{271211719}{73513440}$ \\ [30pt]\hline

8 & ${_3}F_{2}\left(\begin{array}{r} 1, 1, 10;\\
~\\
2, 11 ;\end{array}\  1\right) =\frac{7129}{2268}$ & 18 & ${_3}F_{2}\left(\begin{array}{r} 1, 1, 20;\\
~\\
2, 21 ;\end{array}\  1\right) = \frac{275295799}{73717644}$  \\ [30pt]\hline

9 & ${_3}F_{2}\left(\begin{array}{r} 1, 1,11;\\
~\\
2, 12 ;\end{array}\  1\right) = \frac{81191}{25200}$  &19 & ${_3}F_{2}\left(\begin{array}{r} 1, 1, 21;\\
~\\
2, 22 ;\end{array}\  1\right)= \frac{11167027}{2956096}$  \\ [30pt]\hline

10 &  ${_3}F_{2}\left(\begin{array}{r} 1, 1, 12;\\
~\\
2, 13 ;\end{array}\  1\right) = \frac{83711}{25410}$  &20 & ${_3}F_{2}\left(\begin{array}{r} 1, 1, 22;\\
~\\
2, 23 ;\end{array}\  1\right) =\frac{18858053}{4938024}$  \\ [30pt]\hline

\end{tabular}

\end{table}

\newpage

\begin{table}[h!]
\centering
\caption{( Clausen's summation theorem )}
\vskip.2cm
\label{table-1.2}
\begin{tabular}{|c  |c  |c  |c|c|c|c|c|c|}
\hline S.No. & Clasuen's series & S.No. & Clasuen's series \\ [5pt]\hline
21 & ${_3}F_{2}\left(\begin{array}{r} 1, 1, 23;\\
~\\
2, 24 ;\end{array}\  1\right) =\frac{439143531}{113809696}$ & 31 & ${_3}F_{2}\left(\begin{array}{r} 1, 1, 33;\\
~\\
2, 34 ;\end{array}\  1\right) = \frac{586061125622639}{140027687654400}$     \\  [30pt]\hline

22 & ${_3}F_{2}\left(\begin{array}{r} 1, 1, 24;\\
~\\
2, 25 ;\end{array}\  1\right) =\frac{1332950097}{342075734}$ &32 & ${_3}F_{2}\left(\begin{array}{r} 1, 1, 34;\\
~\\
2, 35;\end{array}\  1\right) = \frac{53676090078349}{12741489961200}$     \\ [30pt]\hline

23 &  ${_3}F_{2}\left(\begin{array}{r} 1, 1, 25;\\
~\\
2, 26 ;\end{array}\  1\right) = \frac{33695573875}{8561966208}$ &33 & ${_3}F_{2}\left(\begin{array}{r} 1, 1, 35;\\
~\\
2, 36 ;\end{array}\  1\right) = \frac{54062195834749}{12752521554240}$     \\ [30pt]\hline

24 & ${_3}F_{2}\left(\begin{array}{r} 1, 1, 26;\\
~\\
2, 27 ;\end{array}\  1\right) = \frac{34052522467}{8580495000}$ &34 & ${_3}F_{2}\left(\begin{array}{r} 1, 1, 36;\\
~\\
2, 37 ;\end{array}\  1\right) = \frac{54437269998109}{12762940281000}$   \\ [30pt]\hline

25 &  ${_3}F_{2}\left(\begin{array}{r} 1, 1, 27;\\
~\\
2, 28 ;\end{array}\  1\right) = \frac{309561680403}{77338861600}$  &35 & ${_3}F_{2}\left(\begin{array}{r} 1, 1, 37;\\
~\\
2, 38 ;\end{array}\  1\right) = \frac{2027671241084233}{472593445833600}$  \\ [30pt]\hline

26 &  ${_3}F_{2}\left(\begin{array}{r} 1, 1, 28;\\
~\\
2, 29 ;\end{array}\  1\right) = \frac{312536252003}{77445096300}$  &36 & ${_3}F_{2}\left(\begin{array}{r} 1, 1, 38;\\
~\\
2, 39 ;\end{array}\  1\right) = \frac{2040798836801833}{472938908878800}$  \\ [30pt]\hline

27 &  ${_3}F_{2}\left(\begin{array}{r} 1, 1, 29;\\
~\\
2, 30 ;\end{array}\  1\right) = \frac{9146733078187}{2248776129600}$  &37 & ${_3}F_{2}\left(\begin{array}{r} 1, 1, 39;\\
~\\
2, 40 ;\end{array}\  1\right) = \frac{2053580969474233}{473266655870400}$ \\ [30pt]\hline

28 & ${_3}F_{2}\left(\begin{array}{r} 1, 1, 30;\\
~\\
2, 31 ;\end{array}\  1\right) =\frac{9227046511387}{2251453244040}$ & 38 & ${_3}F_{2}\left(\begin{array}{r} 1, 1, 40;\\
~\\
2, 41 ;\end{array}\  1\right) = \frac{2066035355155033}{473578015512420}$  \\ [30pt]\hline

29 & ${_3}F_{2}\left(\begin{array}{r} 1, 1,31;\\
~\\
2, 32 ;\end{array}\  1\right) = \frac{288445167734557}{69872686884000}$  &39 & ${_3}F_{2}\left(\begin{array}{r} 1, 1, 41;\\
~\\
2, 42 ;\end{array}\  1\right)= \frac{85205313628946333}{19428841662048000}$  \\ [30pt]\hline

30 &  ${_3}F_{2}\left(\begin{array}{r} 1, 1, 32;\\
~\\
2, 33 ;\end{array}\  1\right) = \frac{581548514594714}{139890941865675}$  &40 & ${_3}F_{2}\left(\begin{array}{r} 1, 1, 42;\\
~\\
2, 43 ;\end{array}\  1\right) =\frac{85691034670497533}{19440406448751600}$  \\ [30pt]\hline

\end{tabular}

\end{table}

\newpage

\begin{table}[h!]
\centering
\caption{( Clausen's summation theorem )}
\vskip.2cm
\label{table-1.3}
\begin{tabular}{|c  |c  |c  |c|c|c|c|c|c|}
\hline S. No. & Clasuen's series  \\ [5pt]\hline
41 & ${_3}F_{2}\left(\begin{array}{r} 1, 1, 43;\\
~\\
2, 44 ;\end{array}\  1\right) =\frac{75614351220200831}{17069625174513600}$         \\  [30pt]\hline

42 & ${_3}F_{2}\left(\begin{array}{r} 1, 1, 44;\\
~\\
2, 45 ;\end{array}\  1\right) =\frac{5853599356775405587}{1315072372819818600}$     \\ [30pt]\hline

43 &  ${_3}F_{2}\left(\begin{array}{r} 1, 1, 45;\\
~\\
2, 46 ;\end{array}\  1\right) = \frac{5884182435213075787}{1315751996785100160}$    \\ [30pt]\hline

44 & ${_3}F_{2}\left(\begin{array}{r} 1, 1, 46;\\
~\\
2, 47 ;\end{array}\  1\right) = \frac{5914085889685464427}{1316402071882326000}$     \\ [30pt]\hline

45 &  ${_3}F_{2}\left(\begin{array}{r} 1, 1, 47;\\
~\\
2, 48 ;\end{array}\  1\right) = \frac{279336945645849479669}{61900150757844484800}$  \\ [30pt]\hline

46 &  ${_3}F_{2}\left(\begin{array}{r} 1, 1, 48;\\
~\\
2, 49 ;\end{array}\  1\right) = \frac{280682601097106968469}{61928185246412349150}$ \\ [30pt]\hline

47 &  ${_3}F_{2}\left(\begin{array}{r} 1, 1, 49;\\
~\\
2, 50 ;\end{array}\  1\right) = \frac{13818010880930033053031}{3035798698036894732800}$ \\ [30pt]\hline

48 & ${_3}F_{2}\left(\begin{array}{r} 1, 1, 50;\\
~\\
2, 51 ;\end{array}\  1\right) =\frac{13881256687139135026631}{3037063614161076772272}$ \\ [30pt]\hline

49 & ${_3}F_{2}\left(\begin{array}{r} 1, 1, 51;\\
~\\
2, 52 ;\end{array}\  1\right) =\frac{13943237577224054960759}{3038278925731369320000}$ \\ [30pt]\hline

50 & ${_3}F_{2}\left(\begin{array}{r} 1, 1, 52;\\
~\\
2, 53 ;\end{array}\  1\right) =\frac{14004003155738682347159}{3039447494548958308200}$ \\ [30pt]\hline

\end{tabular}

\end{table}

\newpage
\section{Application in digamma function}

\begin{table}[h!]
\centering
\caption{ Digamma function for positive integers }
\label{table-1.4}
\begin{tabular}{|c  |c  |c  |c|c|c|c|c|c|}
 \hline
z & $\psi(z)= \frac{\Gamma^{'}{(z)}}{\Gamma{(z)}}$ & z & $\psi(z)= \frac{\Gamma^{'}{(z)}}{\Gamma{(z)}}$ \\  [5pt]\hline

 1 & $-\gamma$  &27 & $ -\gamma+\frac{34395742267}{8923714800}$   \\  [5pt]\hline

 2 & $-\gamma + 1$  &28&  $-\gamma+\frac{312536252003}{80313433200}$     \\  [5pt]\hline

 3 & $-\gamma + \frac{3}{2}$&29 &  $-\gamma+\frac{315404588903}{80313433200}$   \\ [5pt]\hline

 4  & $-\gamma + \frac{11}{6}$ &30&   $-\gamma+\frac{9227046511387}{2329089562800}$       \\  [5pt]\hline

 5  & $-\gamma + \frac{25}{12}$ &31&     $-\gamma+\frac{9304682830147}{2329089562800}$   \\ [5pt]\hline

 6  & $-\gamma + \frac{137}{60}$ &32&    $-\gamma+\frac{290774257297357}{72201776446800}$  \\  [5pt]\hline

 7  & $-\gamma + \frac{49}{20}$ &33&   $-\gamma+\frac{586061125622639}{144403552893600}$      \\  [5pt]\hline

 8  & $-\gamma + \frac{363}{140}$ &34&    $-\gamma+\frac{53676090078349}{13127595717600}$    \\  [5pt]\hline

 9  & $-\gamma + \frac{761}{280}$ &35&   $-\gamma+\frac{54062195834749}{13127595717600}$     \\  [5pt]\hline

 10 & $-\gamma + \frac{7129}{2520}$ &36&    $-\gamma+\frac{54437269998109}{13127595717600}$   \\  [5pt]\hline

 11 & $-\gamma + \frac{7381}{2520}$ &37&   $-\gamma+\frac{54801925434709}{13127595717600}$  \\ [5pt]\hline

 12 & $-\gamma + \frac{83711}{27720}$  &38&  $-\gamma+\frac{2040798836801833}{485721041551200}$    \\ [5pt]\hline

 13 & $-\gamma + \frac{86021}{27720}$  &39&      $-\gamma+\frac{2053580969474233}{485721041551200}$  \\ [10pt]\hline

 14 & $-\gamma + \frac{1145993}{360360}$ &40&$-\gamma+\frac{2066035355155033}{485721041551200}$    \\ [10pt]\hline

 15 & $-\gamma + \frac{1171733}{360360}$ &41& $-\gamma+\frac{2078178381193813}{485721041551200}$ \\ [10pt]\hline

 16 & $-\gamma + \frac{1195757}{360360}$ &42& $-\gamma+\frac{85691034670497533}{19914562703599200}$ \\ [10pt]\hline

 17 & $-\gamma + \frac{2436559}{720720}$ &43&  $-\gamma+\frac{12309312989335019}{2844937529085600}$ \\ [10pt]\hline

 18 & $-\gamma + \frac{42142223}{12252240}$  &44&  $-\gamma+\frac{532145396070491417}{122332313750680800}$\\ [10pt]\hline

 19 & $-\gamma + \frac{14274301}{4084080}$  &45&  $-\gamma+\frac{5884182435213075787}{1345655451257488800}$\\ [10pt]\hline

 20 & $-\gamma + \frac{275295799}{77597520}$ &46&  $-\gamma+\frac{5914085889685464427}{1345655451257488800}$ \\ [10pt]\hline

 21 & $-\gamma + \frac{55835135}{15519504}$ &47&  $-\gamma+\frac{5943339269060627227}{1345655451257488800}$ \\ [10pt]\hline

 22 & $-\gamma +\frac{18858053}{5173168}$ &48&   $-\gamma+\frac{280682601097106968469}{63245806209101973600}$ \\ [10pt]\hline

 23 & $-\gamma + \frac{19093197}{5173168}$ &49&  $-\gamma+\frac{282000222059796592919}{63245806209101973600}$ \\ [10pt]\hline

 24 & $-\gamma + \frac{444316699}{118982864}$ &50&  $-\gamma+\frac{13881256687139135026631}{3099044504245996706400}$  \\ [10pt]\hline

25 & $-\gamma + \frac{1347822955}{356948592}$ &51& $-\gamma+\frac{13943237577224054960759}{3099044504245996706400}$ \\ [10pt]\hline

 26 & $-\gamma+\frac{34052522467}{8923714800}$ &52& $-\gamma+\frac{14004003155738682347159}{3099044504245996706400}$ \\ [10pt]\hline

\end{tabular}

\end{table}

\end{document}